\documentclass[11pt,twoside]{article}
\usepackage{latexsym}
\usepackage{amssymb,amsbsy,amsmath,amsfonts,amssymb,amscd}
\usepackage{epsfig, graphicx}
\usepackage{hyperref}
\usepackage{color}
\newcommand{\commentout}[1]{}
\newcommand{\R}{\mathbb{R}}
\newcommand{\N}{\mathbb{N}}

\newcommand {\al} {\alpha}

\newcommand {\sg} {\sigma}

\newcommand {\Chi} {{\bf \raise 2pt \hbox{$\chi$}} }

\newcommand {\dv}  { {\rm div} }
\newcommand {\f}   {\frac}
\newcommand {\p}   {\partial}
\newcommand{\dis}{\displaystyle}
\newcommand {\proof} {\noindent {\bf Proof}. }
\newcommand{\beq}{\begin{equation}}
\newcommand{\eeq}{\end{equation}}
\newcommand{\bea} {\begin{array}{rl}}
\newcommand{\eea} {\end{array}}
\newcommand{\bepa}{\left\{ \begin{array}{l}}
\newcommand{\eepa} {\end{array}\right.}
\newtheorem{theorem}{Theorem}[section]
\newtheorem{lemma}[theorem]{Lemma}

\newcommand{\qed}{{ \hfill
                       {\unskip\kern 6pt\penalty 500 \raise -2pt\hbox{\vrule\vbox to 6pt{\hrule width 6pt
                       \vfill\hrule}\vrule} \par}   }}
\title{\Large \bf Incompressible limit of mechanical model of tumor growth with viscosity}

\author{
Beno\^ \i t Perthame\thanks{Sorbonne Universit\'es, UPMC Univ Paris 06, UMR 7598, Laboratoire Jacques-Louis Lions, F-75005, Paris, France} 
\thanks{CNRS, UMR 7598, Laboratoire Jacques-Louis Lions, F-75005, Paris, France}
\thanks{INRIA-Paris-Rocquencourt, EPC MAMBA, Domaine de Voluceau, BP105, 78153 Le Chesnay Cedex, France} 
\thanks{Emails:~benoit.perthame@upmc.fr, nicolas.vauchelet@upmc.fr}
\and
Nicolas Vauchelet\footnotemark[1]  \footnotemark[2]  \footnotemark[3] \footnotemark[4]
}

\date{\today}

\begin{document}
\maketitle
\pagestyle{plain}
\begin{abstract} 
Various models of tumor growth are available in the litterature. A first class describes the evolution of the cell number density when considered as a continuous visco-elastic material with growth. A second class, describes the tumor as a set and rules for the free boundary are given related to the classical Hele-Shaw model of fluid dynamics. 
\\

Following the lines of previous papers where the material is described by a purely elastic material, or when active cell motion is included, we make the link between the two levels of description considering the `stiff pressure law' limit. Even though viscosity is a regularizing effect, new mathematical difficulties arise in the visco-elastic case  because estimates on the pressure field are weaker and do not imply  immediately compactness.  For instance, traveling wave solutions and numerical simulations show
that the pressure may be discontinous in space which is not the case for the elastic case.
\end{abstract}

\\
\noindent {\bf Mathematical Classification numbers:} 35J60; 35K57; 74J30; 92C10;
\pagenumbering{arabic}

\section{The cell model with visco-elastic flow}

We consider a mechanical model of tumor growth considered as a visco-elastic media. 
 We denote the number density of tumor cells by $n(x, t)$,
the pressure by $p(x, t)$ and we assume a Brinkman flow that means the macroscopic velocity field is given by $\nabla W$ for a potential $W$ closely related to the pressure. With these assumptions, the model for tumor growth writes
\begin{eqnarray}
&\dis \p_t n_k - \dv (n_k \nabla W_k) = n_k G\big(p_k\big), \label{eqvis:n}
\\[2mm]
&\dis -\nu \Delta W_k + W_k = p_k(x,t):= \Pi_k(n_k), \label{eqvis:W}
\end{eqnarray}
where we choose the pressure law given by:
\beq\label{pneg}
\Pi_k(n) = \f{k}{k-1} n^{k-1}, \qquad k>2.
\eeq
Following \cite{byrne-drasdo,JJP}, we assume that growth is directly related to the 
pressure through a function $G(\cdot)$ that satisfies 
\beq \label{hypG}
G\in C^1(\R), \qquad G'(\cdot)\leq -\alpha <0, \qquad  G(P_M) =0 \quad \text{ for some }\; P_ M>0.
\eeq
The pressure $P_M$ is usually called the {\it homeostatic pressure} . 
We complete equation~\eqref{eqvis:n}, \eqref{eqvis:W}  with a family of initial data $n^0_k$ satisfying (for some constant $C$ independent of $k$)
\beq \label{id1}
0 \leq n^0_k, \qquad \Pi_k(n^0_k ) \leq  p_M, \qquad \| n^0_k\|_{L^1(\R^d)} \leq C.
\eeq

The viscosity coefficient, $\nu >0$, is supposed to be constant; when viscosity is neglected, that means equation \eqref{eqvis:W} with $\nu =0$,  we recover Darcy's law for which an important literature is available, see \cite{CBCB,SC, friedman, RCM, CLLW, Lowengrub_survey, PQV, PQTV} and the references therein. In that case only friction with the cell surrounding (extra-cellular matrix) is considered. Viscosity is a way to represent friction between cells themselves, considered as a Newtonian fluid and Brinkman's law has been derived rigorously for inhomogenous materials \cite{allaire}. Viscoelastic models for tumor growth, based on Stokes' or Brinkman's law have also been used in the context of tumor growth is  \cite{ZWC, BCGRS,JJP} with a major difference, namely the pressure does not follow a law-of-state \eqref{hypG} but follows from the tissue incompressibility.  However, Stokes' or Brinkman's law are also used considering the tissue as `compressible' \cite{BiJu, Ba_ju}. To use Laplacian in \eqref{eqvis:W}, rather than Stokes viscosity terms, is to simplify the presentation and presentation of the mathematical ideas.   Indeed, this is not central for our  aim here, which  is to explain the derivation of such `incompressible' models from the `compressible' equations. Note that the theory of mixtures allows for a general formalism containing both Darcy's law and Brinkman's law \cite{byrneking,AmPr,preziosi_tosin}. 
\\

Our interest is in the  `stiff pressure law' limit of this model towards a free boundary 
model which generalizes the classical Hele-Shaw equation. That is the limit $k\to \infty$ and we first explain formally what can be expected. The limit uses strongly the equation satisfied by the pressure. Multiplying equation \eqref{eqvis:n} by $\Pi_k'(n)$ and using the chain rule, we deduce
$$
\p_t p_k -n_k \Pi_k'(n_k)\Delta W_k - \nabla p_k \cdot \nabla W_k = n_k\Pi_k'(n_k) G(p_k).
$$
From our choice for the law of state  \eqref{pneg}, we deduce that 
$$
n\Pi_k'(n) = k n^{k-1}= (k-1)\Pi_k(n).
$$
Injecting this expression into the above equation, we deduce that
\beq\label{eq:p}
\begin{cases}
\p_t p_k- \nabla p_k \cdot \nabla W_k
= \f{k-1}{\nu}  p_k Q_k,
\\[5pt]
\dis Q_k = W_k-p_k+ \nu G(p_k):= W_k - H^{-1}(p_k)
\end{cases}
\eeq
where we have defined the function $H$, coming with some properties, as
\beq\label{defH}
H := (I-\nu G)^{-1}, \qquad p_m := H(0)>0; \qquad H \mbox{ is increasing}, \qquad H'(\cdot) <1.
\eeq
Indeed, $G$ is non-increasing and thus $(I-\nu G)$ is invertible on
$[0,P_M]$ onto $[-\nu G(0),P_M]$. Furthermore, notice that $(I-\nu G)' > 1$.
\\

Back to the limit $k\to +\infty$, at least when $p_k$ converges strongly, from \eqref{pneg},  we first find the relation
\beq\label{pliml}
p_\infty (1-n_\infty) = 0.
\eeq
Letting $k\to+\infty$ and asuming we can pass into the limit in all terms, we formally deduce
$$
p_\infty \big( \Delta W_\infty + G(p_\infty) \big) = 0.
$$

Therefore, at the limit we can distinguish between two different regions.
The first region is defined by the set
\beq \label{def:omega}
\Omega(t):=\{p_\infty(\cdot, t)>0\}
\eeq
on which we have the system~:
\begin{eqnarray}
&\dis n_\infty =1, \label{omega:nlim}\\[2mm]
&\dis -\nu\Delta W_\infty + W_\infty = p_\infty, \label{omega:Wlim}\\[2mm]
&\dis \Delta W_\infty + G(p_\infty) = 0. \label{omega:plim}
\end{eqnarray}
Thus the latter system reduces to~:
\begin{eqnarray*}
& n_\infty = 1, \qquad p_\infty = H(W_\infty), \qquad x \in \Omega(t),
\\[2mm]
& -\nu \Delta W_\infty + W_\infty - H(W_\infty) = 0.
\end{eqnarray*}

On the second region, $\R^d\setminus \Omega(t)$, the limiting system writes
\begin{eqnarray*}
&\dis p_\infty = 0, \\[2mm]
&\dis \p_t n_\infty - \dv (n_\infty\nabla W_\infty) = n_\infty G(0), \\[2mm]
&\dis -\nu\Delta W_\infty + W_\infty = 0.
\end{eqnarray*}
\\

To establish rigorously this limit, we need some additional assumption on the initial data. Namely, we need that the family $n^0_k$ is  `well-prepared'. By this, we mean that,  for some open set $\Omega^0$, 
\beq  
 \Pi_k(n^0_k) \underset{ k \to \infty}{\longrightarrow}  p^0_\infty = H(W_\infty) \quad  \text{a.e. in  } \;  \Omega^0, \qquad n^0_k =0  \quad \text{in  } \R^d \backslash \Omega^0 . 
 \label{id2}
\eeq
Note that, with the notation in~\eqref{eq:p}, this assumption implies that $Q^0_k \equiv 0$ and $n^0_k \underset{k\to +\infty}{\longrightarrow}  \mathbf{1}_{\Omega^0}$. 
For this purpose, the latter assumption can be slightly relaxed to $n^0_k  \ll e^{-A/k}$ for all $A>0$ in $ \R^d \backslash \Omega^0$. With our present proof, we need to avoid the existence of a domain where $n^0_k$ remains strictly between $0$ and $1$, a case which we leave open at this stage. 

Our goal is to prove the
\begin{theorem}\label{th:lim}  
Under assumptions \eqref{hypG}, \eqref{id1} and \eqref{id2}, consider a solution of the system \eqref{eqvis:n}-- \eqref{pneg}. After extraction of subsequences, both the density $n_k$ and the pressure $p_k$ converge strongly in $L^1_{loc}\big( (0,T)\times \R^d \big)$, for all $T>0$, as $k\to +\infty$ towards 
respectively $n_\infty$ and $p_\infty$ belonging to $L^1\cap L^\infty\big((0,T) \times \R^d\big)$;
up to a subsequence, $W_k$ converges strongly in $L^1\big((0,T),W^{1,q}_{loc}(\R^d) \big)$, 
for all $q\geq 1$, towards $W_\infty$. Moreover, these functions satisfy
\beq\label{eq:nlim}
\p_t n_\infty - \dv (n_\infty \nabla W_\infty) = n_\infty G(p_\infty),\qquad
n_\infty(t=0) = n^0_\infty=  \mathbf{1}_{\{ \Omega^0 \}} ,
\eeq
\beq\label{eq:Wlim}
-\nu \Delta W_\infty + W_\infty = p_\infty,
\eeq
\beq\label{eq:plim}
p_\infty=H(W_\infty) \mathbf{1}_{\{p_\infty>0\}}, \quad \qquad p_\infty (1-n_\infty) = 0,
\eeq
\beq\label{relcomp}
p_\infty \big( p_\infty -W_\infty - \nu G(p_\infty) \big) = 0, \quad \mbox{ a.e.}
\eeq
\end{theorem}

The first relation in~\eqref{eq:plim} is equivalent to the statement~\eqref{relcomp} and replaces the usual  `complementary relation' in Hele-Shaw flow, $p_\infty (\Delta p_\infty + G(p_\infty))$, see \cite{PQV, PQTV, Elliot_j}.

Because the function $H(\cdot)$ does not vanish, we conclude from the first relation in~\eqref{eq:plim},  that $p_\infty$ is discontinuous. This is a major difference with elastic materials (Darcy's law), then $p_\infty$ is continous in space, and this is illustrated by traveling wave solutions we build in Section~\ref{sec:tw}. The pressure jump is however related to the potential $W_\infty$, a difference with models including surface tension where the jump is related to the free boundary curvature, see \cite{ABC, ES} and the reference therein.
\\

We first prove Theorem~\ref{th:lim} in several steps. In a first step, we derive a priori estimates. Because they do not give compactness for the pressure, we analyze possible oscillations using a kinetic formulation. From properties of solutions of the corresponding kinetic equation, we conclude that strong compactness occurs. All these steps are in Section~\ref{sec:proof}. The one dimensional traveling wave profiles are presented in Section \ref{sec:tw} with numerical illustrations. The final Section is devoted to a conclusion and presentation of some perspective.

\section{Proof of the Hele-Shaw limit}
\label{sec:proof}

We divide the proof of our main result Theorem~\ref{th:lim} in several steps. We begin with several bounds which are useful for the sequel. Then, in order to prove strong convergence of the pressure $p_k$, we analyze possible oscillations using the kinetic formulation of \eqref{eq:p} in the spirit of \cite{BPALD}.

\subsection{Estimates}

\begin{lemma}[A priori estimates]\label{lem:estim}
Under previous assumptions, for all $T>0$, the uniform bounds with respect to $k$ hold
$$
\begin{array}{l}
\dis n_k, \; p_k \text{ and } W_k \in  L^\infty \big((0,T); L^1\cap L^\infty( \R^d) \big),  \qquad p_k \leq P_M, 
 \\[2mm]
\dis W_k \in L^\infty \big((0,T); W^{1,q}(\R^d)\big), \mbox{ for } 1 \leq q \leq \infty, \qquad \dis D^2 W_k \in L^\infty \big((0,T); L^{q}(\R^d)\big), \mbox{ for } 1 < q < \infty, 
\\[2mm]
\p_t W_k \in L^1\big((0,T); L^{q}(\R^d)\big), \mbox{ for } 1 \leq q \leq \infty, \quad  \p_t \nabla W_k \in L^1\big((0,T); L^{q}(\R^d)\big), \quad  \mbox{ for } 1< q <\f{d}{d-1}.
\end{array}
$$
For some nonnegative constant $C(T)$, independent of $k$, we have
\beq\label{estim1}
k \int_0^T \int_{\R^d}  p_k\; \big| p_k-W_k -  \nu G(p_k)\big|\, dxdt \leq C(T).
\eeq
\end{lemma}

We can draw several consequences of this Lemma. First, after extracting subsequences, it is immediate that the following convergences hold as $k\to \infty$: 
$$
n_k \to n_\infty \leq 1, \qquad p_k \to p_\infty \leq P_M \quad \text{weakly}-\star \mbox{ in } L^\infty \big((0,+\infty) \times \R^d \big),
$$
and these limits belong to  $L^\infty\big((0,T); L^1(\R^d)\big)$ for all $t>0$. Also, we have
$$
W_k \to  W_\infty < P_M , \quad \nabla W_k \to \nabla W_\infty    \quad \text{locally in } L^q  \big((0,T) \times \R^d \big),  \quad 1 \leq q < \infty .
$$ 
Passing to the limit 
in \eqref{eqvis:W} and in the left hand side of \eqref{eqvis:n}, we get
\beq\label{eqWinf}
-\nu \Delta W_\infty + W_\infty = p_\infty.
\eeq

The second consequence concerns the backward flow with velocity $\nabla W_k$ defined as 
\beq
\f{d}{ds} X^{(k)}_{(x,t)}(s)= - \nabla W_k(X^{(k)}_{(x,t)}(s),s) , \qquad X^{(k)}_{(x,t)}(t)= x,
\label{flowk}
\eeq
as well as the forward flow
\beq
\f{d}{dt} Y^{(k)}_{(x)}(t)= - \nabla W_k(Y^{(k)}_{(x)}(t),t) , \qquad Y^{(k)}_{(x)}(t=0)= x .
\label{fflowk}
\eeq
Even though, $ \nabla W_k$ is not uniformly Lipschitz continuous but slightly less, and  according to DiPerna-Lions theory \cite{DiPernaLions}, these flow are well defined a.e. and, after extraction of subsequences as in Lemma~\ref{lem:estim}, it converges a.e. to the limiting flows defined by~\eqref{flowlimit} for the backward flow and by \eqref{flowlimitY} for the forward flow.
\\

The third conclusion uses a combination of the above flow with equation~\eqref{eq:p}. We have
\beq
p_k(x,t)=0 \quad \text{for } x \in \R^d \backslash \Omega^k(t), \qquad \Omega^k(t) = Y^{(k)}(t)[\Omega^0]
\label{compactsupp}
\eeq

\proof \\
{\bf 1st step. A priori bounds in $L^1\cap L^\infty$.} Clearly $n_k$ is nonnegative provided $n_k(t=0)\geq 0$. Integrating, we deduce
a bound for $n_k$ in  $L^\infty \big((0,T); L^1 ( \R^d) \big)$, uniformly with respect to $k$.
\\

By definition of $p_k$ in  \eqref{pneg}, we clearly have that $\Pi_k'(n_k)\geq 0$ when $k>1$.
We can apply the maximum principle of \cite[Lemma 2.1]{TVCVDP} to obtain the uniform bound 
$$
0\leq p_k \leq P_M.
$$ 
Therefore, still using relation~\eqref{pneg}, we have $n_k=\big(\f{k-1}{k} p_k\big)^{1/(k-1)}$
and $n_k$ is uniformly bounded in $L^\infty \big( (0,+\infty ) \times \R^d \big)$.
Then, writing $p_k\leq n_k \|n_k\|_\infty^{k-2}$, we deduce an uniform bound of 
$(p_k)_k$ in $L^\infty \big((0,T); L^1 ( \R^d) \big)$.
\\ \\
{\bf 2nd step. Representation of $W_k$.} 
Using elliptic regularity on \eqref{eqvis:W}, we conclude that for all $t\in [0,T]$, $W_k(t,\cdot)$ is
bounded in $W^{2,q}(\R^d)$. Moreover, denoting by $K$ the fondamental solution
of $-\nu\Delta K+K=\delta_0$, we have 
\beq\label{def:K}
W_k=K\star p_k , \qquad K(x)= \f{1}{4\pi} \int_0^{\infty} e^{-\big(\pi \f{|x|^2}{4s\nu} + \f{s}{4\pi}\big)}
\,\frac{ds}{s^{d/2}}.
\eeq
We recall that 
$$ \begin{cases}
K\in L^{q}(\R^d), \qquad \forall 1\leq q < \frac{d}{d-2}, \qquad  (1\leq q \leq +\infty\; \text{for } d=1),
\\[5pt]
\nabla K \in L^{q}(\R^d), \quad \forall 1\leq q<\f{d}{d-1},
\end{cases}
$$
and that  $K \geq 0$, $\int_{\R^d} K(x)\,dx =1$,  which we use below. 
\\

Taking the convolution of \eqref{eq:p}, we deduce
\beq\label{eqt:W}
\p_t W_k = K\star [\nabla p_k\cdot \nabla W_k + \f{k-1}{\nu} p_k Q_k].
\eeq
\\
{\bf 3rd step. Bounds on $Q_k$.} 
Then, by definition of $Q_k$ and using \eqref{eq:p}, we compute
$$
\p_t Q_k - \nabla Q_k\cdot\nabla W_k  + \f{k-1}{\nu} \big( 1 -\nu G'(p_k)\big) p_k Q_k 
= -|\nabla W_k|^2 + K\star [\nabla p_k\cdot \nabla W_k + \f{k-1}{\nu} p_k Q_k].
$$
Therefore, from a standard computation, we deduce
$$
\begin{array}{l}
\dis \p_t |Q_k| - \nabla |Q_k| . \nabla W_k  + \f{k-1}{\nu} \Big(1 -\nu G'(p_k)\Big) p_k\; |Q_k|  
\\[2mm]
\dis \qquad \leq |\nabla W_k|^2 +|K\star [\nabla p_k . \nabla W_k] | + \f{k-1}{\nu} K\star [p_k\; |Q_k|].
\end{array}
$$
We may integrate in $x$ and $t$. Because $p_k$ and $W_k$ are uniformly bounded in $L^1\cap L^\infty$, 
and $|G'|\geq \al$ from \eqref{hypG}, we find
$$\begin{array}{rl}
  \dis \al (k-1) \dis  \int_0^T \int_{\R^d}   p_k |Q_k|  \,dxdt &\leq \dis \int_{\R^d} |Q_k(x,0)|dx  -  \int_{\R^d} |Q_k(x,T)|dx + \dis  \int_0^T \int_{\R^d} |\nabla W_k|^2\; dx dt
\\[2mm]
& +  \dis \int_0^T \int_{\R^d}  \big[- |Q_k|  \; \Delta W_k   +|K\star [\nabla p_k \cdot \nabla W_k] |\big]\,dxdt.
\end{array}
$$
The three first terms in the right hand side are all controlled uniformly and, to conclude the bound \eqref{estim1}, we have to estimate last two terms. Using \eqref{eqvis:W}, the first term is 
$$
-\int_0^T\int_{\R^d}  |Q_k| \Delta W_k \,dxdt=  \f 1 \nu  \int_0^T\int_{\R^d}  |Q_k|\;  ( p_k -W_k)\,dxdt \leq  \f 1 \nu  \int_0^T\int_{\R^d}  |Q_k|\; p_k \,dxdt ,
$$
and this term is controlled, for $k$ large enough, by the $\al k$ term in the left hand side.
The second term is 
$$
K\star [\nabla p_k \cdot \nabla W_k] =  \nabla K\star [p_k\cdot \nabla W_k] - K \star [ p_k \Delta W_k].
$$
Using the uniform bounds on $p_k$, we have that
$p_k \cdot \nabla K* p_k$ is uniformly bounded, with respect to $k$, in $L^\infty \big( (0,T); L^q(\R^d) \big)$, $1 \leq q \leq \infty$, 
and thus, $\nabla K\star [p \cdot \nabla W_k]$ is also uniformly bounded in $L^\infty \big( (0,T); L^q(\R^d) \big)$, $1 \leq q \leq \infty$. Finally, 
$p_k \Delta W_k$ is also uniformly bounded in $L^\infty \big( (0,T); L^q(\R^d) \big)$, $1 \leq q \leq \infty$.
This immediately concludes the proof of estimates \eqref{estim1}. 
\\
\\
{\bf 4th step. Estimate on $\p_t W_k$.} 
Finally, using the above estimate  and equation \eqref{eqt:W}, we deduce that $\p_t W_k$ is 
uniformly bounded with respect to $k$ in $L^\infty \big( (0,T); L^q(\R^d) \big)$, $1 \leq q \leq \infty$.


For the estimate for $\p_t \nabla W_k$, we can use again the above calculation and write 
$$
\p_t \nabla W_k =  - D^2 K\star [p_k\cdot \nabla W_k] + \nabla K \star [ p_k \Delta W_k] + \f{k-1}{\nu}  \nabla K*[ p_k Q_k ].
$$  
Since $D^2 K$ is a bounded operator in $L^1$, we conclude the last bound in Lemma~\ref{lem:estim}.
\qed

\subsection{Which oscillations for the pressure?}

We deduce from Lemma~\ref{lem:estim} that, up to a subsequence, the sequence $(W_k)_k$
converges strongly in $L^1((0,T),W^{1,q}_{loc})$.
However, we only get weak convergence for the pressure $(p_k)_k$ and the 
density $(n_k)_k$.
Here, we give an argument showing that the only obstruction to strong compactness, 
is oscillations of $p_k$ between the values $p_k \approx 0$ and $p_k \approx H(W_\infty)$.

\begin{lemma}\label{lem2}
Let $T>0$ and let $H$ be defined in \eqref{defH} with the assumptions \eqref{hypG}.
Consider real numbers $\beta_1 >0$, $\beta_2 >0$ small enough,  and let $p_k$ be as in Lemma \ref{lem:estim}, then we have
$$
{\rm meas}\big\{ \beta_1\leq p_k(x,t) \leq H(W_\infty(x,t)) - \beta_2 \big\} \underset{k\to +\infty}{\longrightarrow} 0,
$$
$$
{\rm meas}\big\{   p_k(x,t) \geq H(W_\infty(x,t)) + \beta_2 \big\} \underset{k\to +\infty}{\longrightarrow} 0,
$$
where $\rm meas$ denotes the Lebesgue measure.
\end{lemma}

\proof
Let $0<\beta_1<\beta_2<p_m$, $p_m$ being defined in \eqref{defH}, we have for all $k\in \N$
\beq\label{lem2:eq1}
\int_0^T \int_{\R^d} \mathbf{1}_{\{\beta_1 \leq p_k\leq H(W_\infty) - \beta_2 \}}\,dxdt \leq
\int_0^T \int_{\R^d} \f{p_k}{\beta_1}\mathbf{1}_{\{ \beta_1\leq p_k\leq H(W_\infty) - \beta_2 \}}\,dxdt . 
\eeq
From assumption \eqref{hypG}, the function $I-\nu G$ is increasing and by definition \eqref{defH}, 
$(I-\nu G)(H(W_\infty))= W_\infty \geq 0$ (the nonnegativity is because $W_\infty$ is a solution of \eqref{eqWinf}).  
Therefore, on the set $\{p_k\leq  H(W_\infty) - \beta_2 \}$, we have, for some $\omega(\beta_2) >0$,
$$
(I-\nu G)(p_k) \leq (I-\nu G)(H(W_\infty) - \beta_2) \leq W_\infty -\omega(\beta_2), 
$$
$$
W_\infty - (I-\nu G)(p_k)  \geq \omega(\beta_2) .
$$
Thus we can estimate
$$\begin{array}{r}
\dis \int_0^T \int_{\R^d} p_k \mathbf{1}_{\{\beta_1\leq p_k\leq H(W_\infty)-\beta_2\}}\,dxdt  \leq \dis
\f{1}{\omega(\beta_2) } \iint_{\{\beta_1\leq p_k\leq  H(W_\infty) - \beta_2\}} p_k|(I-\nu G)(p_k)-W_\infty| \,dxdt
\\[5pt]
\qquad \qquad \dis \leq 
\f{1}{\omega(\beta_2) } \left[ \int_0^T\int_{\R^d} p_k|(I-\nu G)(p_k)-W_k| \,dxdt + 
\iint_{\{\beta_1\leq p_k\}} p_M |W_\infty - W_k| dx dt \right].
\end{array}
$$

Additionally, using estimate \eqref{estim1}, and the strong convergence of $W_k$, we deduce that
\beq\label{estim2}
\lim_{k\to +\infty} \int_0^T \int_{\R^d} p_k \mathbf{1}_{\{\beta_1\leq p_k\leq H(W_\infty)-\beta_2\}}\,dxdt = 0.
\eeq
We notice, for future use, that in the same spirit we also have that
\beq\label{estim3}
\lim_{k\to +\infty} \int_0^T \int_{\R^d} p_k \mathbf{1}_{\{p_k\leq H(W_\infty)-\beta_2\}}\,dxdt = 0.
\eeq
Thus estimates \eqref{lem2:eq1}--\eqref{estim2} prove the first statement of Lemma~\ref{lem2}.
\\

The second statement can be proved in the same way.
\qed

\subsection{Strong convergence of the pressure}

However, we need strong convergence to recover the asymptotic limit,
in particular the equation satisfied by $p_\infty$.
A difficulty here is that we do not have estimates on the derivatives on $p$, unlike in \cite{PQV,PQTV}.
Then we develop another strategy based on estimate \eqref{estim1} to obtain the 
following strong convergence result~:
\begin{lemma}[Strong convergence of $p_k$]\label{lem:convp}
Up to a subsequence, $p_k$ converges strongly locally  in $L^1 \big((0,T) \times \R^d \big)$
towards $p_\infty$. 
Moreover, $p_\infty=H(W_\infty) \mathbf{1}_{\{p_\infty>0\}}$ a.e.

Furthermore, we have 
$$
\Omega(t)=\{p_\infty(\cdot, t)= H(W_\infty(\cdot, t))\} = \R^d\backslash \{p_\infty(\cdot, t) =0\}
$$ 
is the image of $\Omega^0$ by the limiting flow $Y_{(x)}(t)$, defined by
\beq\label{flowlimitY}
\f{d}{dt} Y_{(x)}(t)= - \nabla W_\infty(Y_{(x)}(t),t) , \qquad Y_{(x)}(t=0)= x.
\eeq

Finally, we have for all $T>0$,
\beq \label{eq:muvanish}
k  \int_0^T \int_{\R^d}p_k (x,t)  |Q_k(x, t)| dx dt  \underset{k\to +\infty}{\longrightarrow} 0. 
\eeq
\end{lemma}
\proof
The strategy is to pass to the limit in the equation~\eqref{eq:p} for $p_k$ and to combine this information with the possible oscillations of $p_k$ as described by Lemma~\ref{lem2}. For that, we need a representation of the weak limit of $p_k$ which we can obtain thanks to a kinetic representation. 
\\
\\
\noindent {\bf 1st step. Representation of nonlinear weak limits.} Our first result is that there is a measurable function $0\leq f(x,t)\leq 1$ such that for all smooth function $S: [0, \infty) \to \R$, we have, up to a subsequence,
\beq\label{lim:Sp}
S(p_{k}) \underset{k\to +\infty}{\rightharpoonup} S(0)(1-f) + S(H(W_\infty)) f,
\eeq
and
\beq\label{lim:Spp}
S(0)(1-f) + S(H(W_\infty)) f = \int_0^\infty S'(\xi) \chi(\xi) \,d\xi + S(0) , \qquad 
\chi(x,\xi,t) = f(x,t) \mathbf{1}_{\{0 < \xi < H(W_\infty(x,t))\}}.
\eeq
Interpreted in terms of Young measures, this means that $p_k$ oscillates between the values $0$ and $H(W_\infty(x,t))$ with the weights $1-f(x,t)$ and $f(x,t)$. Notice that for $S(p)=p$, we find
\beq \label{lim:pWinfty}
p_\infty = f \; H(W_\infty).
\eeq 

To prove these results, we define 
$$
\chi_k (x, \xi, t) =   \mathbf{1}_{\{0 < \xi < p_k(x,t)\}} 
$$
and we write 
\beq\label{eq:Spk}
S(p_k) -S(0) = \int_0^\infty S'(\xi) \chi_k (x, \xi, t)  d \xi.
\eeq
We can extract a subsequence, still denoted $(p_k)_k$, such  that  $\mathbf{1}_{\{0 < \xi < p_{k}\}} $ converges in $L^\infty((0,\infty)\times \R^d)-weak\star$ towards a function $\chi(x,\xi,t) $ which satisfies $0 \leq \chi(x,\xi,t) \leq 1$. 
Then $S(p_{k})$ converges weakly to $S(0)+\int_0^\infty S'(\xi) \chi(x,\xi,t) d\xi $.  

We define, 
$$
f(x,t) = {\rm w\! -\! lim}\;   \mathbf{1}_{ \{ p_{k}(x,t) \geq p_m/2 \} }
$$
where we recall that $p_m$ is defined in \eqref{defH}.
Since $H(W_\infty ) > p_m$, we  may use Lemma~\ref{lem2} to conclude \eqref{lim:Sp}--\eqref{lim:Spp}.

\medskip
\noindent {\bf 2nd step. Equation satisfied by $\chi_k$.}
We use the equation \eqref{eq:p}
$$
\p_t p_k - \nabla p_k\cdot\nabla W_k = \f{k-1}{\nu} p_k Q_k, \qquad
Q_k =  W_k-p_k + \nu G(p_k).
$$
For any  function  $S \in C^2(\R; \R)$, multiplying it by $S'(p_k)$ leads to
$$
\p_t S(p_k) - \nabla S(p_k)\cdot\nabla W_k = (k-1) p_k Q_k S'(p_k).
$$
Denoting $\delta$ the Dirac mass, we can rewrite the later equation as
\beq\label{eq:Sp}
\p_t \int_0^\infty S'(\xi) \chi_k d\xi - \nabla  \int_0^\infty S'(\xi)\chi_k d\xi \cdot\nabla W_k  = \int_0^\infty  S'(\xi) \mu_k(x, \xi, t) \,d\xi,
\eeq
\beq\label{eq:muk}
\mu_k(x, \xi, t) :=  \f{k-1}{\nu} p_kQ_k \delta_{\{\xi =p_k \}} =   \f{k-1}{\nu} p_k [W_k-p_k+ \nu G(p_k)] \delta_{\{\xi =p_k \}}.
\eeq 
Eliminating the test function $S'(\cdot)$, this is equivalent to write 
\beq\label{eq:SkTrueStrong}
\p_t \chi_k -\nabla \chi_k \cdot \nabla W_k =  \mu_k .
\eeq
\\

However, this formula is not enough to pass to the limit $k\to \infty$ and we need the divergence form,
$$
\p_t S(p_k) - \dv\big [ S(p_k)\nabla W_k]+  S(p_k) \f{W_k- p_k}{\nu} = (k-1) p_k Q_k S'(p_k)
=\int_0^\infty S'(\xi) \mu_k(d\xi).
$$
Therefore, using \eqref{eq:Spk} and the fact that $S(p_k)p_k = \int_0^\infty \big(S(\xi)+\xi S'(\xi)\big)\chi_k\,d\xi$, we have
\beq\label{eq:Sk}
\int_0^\infty S'(\xi) \big[ \p_t \chi_k - \dv [\chi_k \cdot \nabla W_k] + \chi_k \f{W_k- \xi}{\nu} \big] \,d\xi - \int_0^\infty \f{S(\xi)-S(0)}{\nu}  \chi_k d\xi=
\int_0^\infty S'(\xi) \mu_k(d\xi).
\eeq

Because $\chi_k(\xi)= - \frac{\p}{\p \xi} \int_\xi^{\infty} \chi_k(x, \eta, t)d\eta$, and integrating by parts, 
we have
$$
 \int_0^\infty \f{S(\xi)-S(0)}{\nu}  \chi_k d\xi =  \int_0^\infty \f{S'(\xi)}{\nu}   \int_\xi^{\infty} \chi_k(x, \eta, t)d\eta d\xi.
$$
Therefore, \eqref{eq:Sk} is equivalent to our final formulation
\beq\label{eq:SkTrue}
\p_t \chi_k -\dv [ \chi_k \nabla W_k] + \chi_k  \f{W_k- \xi}{\nu}  - \f 1 \nu \int_\xi^{\infty} \chi_k(x, \eta, t) d\eta =  \mu_k.
\eeq

One can simplify this relation and write
$$
\p_t \chi_k -\dv [ \chi_k  \nabla W_k] + \chi_k  \f{W_k- \xi}{\nu}  - \f {(p_k- \xi)_+}{ \nu} =  \mu_k.
$$
Finally, \eqref{eq:Sk} is equivalent to
$$
\p_t \chi_k -\dv [ \chi_k   \nabla W_k] + \chi_k  \f{W_k- p_k}{\nu} =  \mu_k .
$$
In particular, integrating in $\xi$ we recover the expected formula
$$
\p_t p_k -\dv [ p_k \nabla W_k] + \f{p_k}{\nu}  [W_k - {p_k}]  = \int \mu_k d\xi. 
$$

\medskip
\noindent {\bf 3rd step. Equation satisfied by $f$.} 
We may pass to the limit in \eqref{eq:SkTrue}. For all $T>0$, the sequence $\mu_k$ is  uniformly bounded in $L^1(\R^d\times \R\times [0,T])$ thanks to estimate \eqref{estim1}.
Thus we can extract a subsequence converging, in the weak sense of measures, towards 
a measure denoted $\mu$ in ${\cal M}_b(\R^d\times \R\times [0,T])$.
Because $Q_k(x,\xi,t) = W_k - \xi + \nu G(\xi)$ is positive for $\xi \leq p_m$, we have 
$$
\mu(x,\xi,t) \geq 0 \qquad  \text{for } \;  \xi \leq p_m.
$$

Therefore passing to the limit $k\to +\infty$ into \eqref{eq:SkTrue}, in the sense of distributions, 
$$
\p_t \chi -\dv [ \chi \cdot \nabla W_\infty] + \chi  \f{W_\infty- \xi}{\nu}  - \f 1 \nu \int_\xi^{\infty} \chi(x, \eta, t) d\eta =  \mu .
$$
This last equation can also be written with \eqref{lim:Spp}
$$
\p_t \chi -\dv [ \chi \cdot \nabla W_\infty] + \chi  \f{W_\infty- \xi}{\nu}  - f(x,t) \f {(H(W_\infty)- \xi )_+} \nu  =  \mu ,
$$
and thus
\beq\label{eq:SlimTrue}
\p_t \chi -\dv [ \chi \cdot \nabla W_\infty] + \chi   \f{W_\infty- H(W_\infty) }{\nu} =  \mu .
\eeq

Using the assumption \eqref{id2}, this equation is complemented with the initial condition 
$$
\chi(x,\xi, t=0) = \mathbf{1}_{\Omega^0} \mathbf{1}_{\{0 < \xi < H(W_\infty(x,t=0))\}}
$$
and 
$$
f(x, t=0)= \mathbf{1}_{ \Omega^0 }:= f^0(x) .
$$

It is useful to keep in mind the equivalent form of this equation, 
$$
\p_t \chi - \nabla  \chi \cdot \nabla W_\infty  + \chi  \f{p_\infty- H(W_\infty)}{\nu}   =  \mu \geq 0.
$$
and thus, using~\eqref{lim:pWinfty},  
\beq\label{eq:SlimW}
\p_t \chi - \nabla  \chi \cdot \nabla W_\infty  + \chi  \; H(W_\infty) \f{f - 1}{\nu}   =  \mu \geq 0.
\eeq

We can also integrate \eqref{eq:SlimTrue} and recover
$$
\p_t p_\infty   -\dv [ p_\infty   \cdot \nabla W_\infty  ] + \f{p_\infty  }{\nu}  [W_\infty   -  H(W_\infty  )]  = \int \mu d\xi  . 
$$
\\
\\
{\bf 4th step. The set $\{ g(x,t)=1 \text{ and }  \xi <p_m\}$.}  It is useful to consider the function
$$
g(x,t) = f^0 \big(X_{(x,t)}(s=0) \big) ,  
$$ 
with the characteristics defined by
\beq
\f{d}{ds} X_{(x,t)}(s)= - \nabla W_\infty(X_{(x,t)}(s),s) , \qquad X_{(x,t)}(t)= x.
\label{flowlimit}
\eeq
This function $g$ is the solution  of the transport equation
$$
\p_t g -\nabla g \cdot \nabla W_\infty = 0 , \qquad g^0=f^0. 
$$

Using \eqref{eq:SlimW} and $0 \leq f \leq 1$,  we find 
\beq\label{eq:f}
\p_t f -\nabla f \cdot \nabla W_\infty = \mu (x,\xi, t) +\chi  \; H(W_\infty) \f{1 - f}{\nu} \geq 0  .
\eeq
From the comparison principle, we conclude that $f(x,t) \geq g(x,t)$ and we conclude that,
\beq\label{eq:mu=0} \begin{cases}
f(x,t) = g(x,t) =1, \qquad \text{ in the set}Ê\quad \{g(x,t) = 1 \},
\\[5pt]
\mu (x, \xi, t) =0 \qquad \text{in the set}Ê\quad \{g(x,t) = 1 \text{ and } \xi < p_m \} .
\end{cases}
\eeq
\\
\\
{\bf 5th step. Strong convergence of $p_k$.}  
Another wording for step 4, is that
$$
\Omega(t)= Y_{(x)}(t)[\Omega^0] = \{ p_\infty(\cdot, t) >0\},
$$
with $Y_{(x)}(t)$ the limiting flow of  $Y_{(x)}^{(k)}(t)$ defined in \eqref{fflowk}.  
Indeed, from \eqref{compactsupp} and the strong convergence of the flow, we infer that
$$
p_\infty(\cdot,t) = 0 \quad \text{ in }    Y_{(x)}(t)[ \R^d \backslash \Omega^0].
$$
Then we have $f(x,t)=\mathbf{1}_{\Omega(t)}=\mathbf{1}_{\{p_\infty(x,t)>0\}}$. We recall that
by definition, $f=w-\lim_{k\to +\infty} \mathbf{1}_{\{p_k\geq p_m/2\}}$.

We show that it implies the strong convergence locally in $L^1((0,T)\times \R^d)$
of $p_k$ towards $H(W_\infty) \mathbf{1}_{\{p_\infty>0\}}$. Let $U$ be an open bounded
subset of $\R^d$, we have
\beq\label{ineq:convp}
\int_0^T\int_U |p_k-H(W_\infty) \mathbf{1}_{\{p_\infty>0\}}|\,dx \leq I_k + I\!I_k + I\!I\!I_k,
\eeq 
with
\begin{eqnarray*}
&&\dis I_k=\int_0^T\int_U \mathbf{1}_{\{p_k\geq p_m/2\}} |p_k - H(W_\infty)|\,dx,
\\ &&\dis 
I\!I_k = \int_0^T\int_U \mathbf{1}_{\{p_k < p_m/2\}}p_k\,dx,
\\ &&\dis 
I\!I\!I_k=\int_0^T\int_U H(W_\infty) \big(\mathbf{1}_{\{p_k \geq p_m/2\}} (1- \mathbf{1}_{\{p_\infty>0\}})
+\mathbf{1}_{\{p_k < p_m/2\}}  \mathbf{1}_{\{p_\infty>0\}}\big) \,dx. \\
\end{eqnarray*}
For the first term $I_k$, we have that
$$
\begin{array}{ll}
\dis I_k & \dis \leq \int_0^T\int_U \mathbf{1}_{\{p_k\geq p_m/2\}} |p_k - H(W_k)|\,dx +
\int_0^T\int_U \mathbf{1}_{\{p_k\geq p_m/2\}} |H(W_k) - H(W_\infty)|\,dx   \\[3mm]
& \dis \leq \frac{2}{p_m} \int_0^T\int_U p_k |p_k - H(W_k)|\,dx + 
C \int_0^T\int_U |W_k - W_\infty|\,dx.
\end{array}
$$
Using estimate \eqref{estim1}, we deduce that the first term of the right hand side 
goes to $0$ as $k\to +\infty$. From the local strong convergence of $W_k$ towards $W_\infty$,
the second term of the right hand side converges to $0$ too. We conclude that 
$\lim_{k\to +\infty} I_k = 0$.
Moreover, it has been proved in Lemma \ref{lem2}, see equation \eqref{estim3}, 
that $\lim_{k\to +\infty} I\!I_k=0$.
For the last term, we have, using the fact that $W_\infty$ is bounded in $L^\infty$,
that for some nonnegative constant $C$,
$$
I\!I\!I_k \leq C \int_U \Big(\mathbf{1}_{\{p_k \geq p_m/2\}} \big(1- \mathbf{1}_{\{p_\infty>0\}}\big) + \big(1-\mathbf{1}_{\{p_k \geq p_m/2\}} \big) \mathbf{1}_{\{p_\infty>0\}}\Big) \,dx
$$
We have shown in the 4th step above that $\mathbf{1}_{\{p_k \geq p_m/2\}}$ converges weakly towards
$\mathbf{1}_{\{p_\infty>0\}}$. Then passing to the limit $k\to +\infty$ in the latter inequality,
we deduce that $\lim_{k\to +\infty} I\!I\!I_k=0$.
We conclude from \eqref{ineq:convp} that, for any open bounded subset $U$,
$$
\int_0^T\int_U |p_k-H(W_\infty) \mathbf{1}_{\{p_\infty>0\}}|\,dx \underset{k\to +\infty}{\longrightarrow} 0.
$$
By uniqueness of the weak limit, we deduce that $p_\infty = H(W_\infty) \mathbf{1}_{\{p_\infty>0\}}$ a.e.
\\
\\
{\bf 6th step. Derivation of \eqref{eq:muvanish}.} From definition \eqref{eq:muk}, this limit is now a consequence of
$$
k  \int_0^T \int_{\R^d} p_k(x,t) |Q_k(x, t)| \,dxdt = \int_0^T \int_{\R^d} \int_{(0,\infty)} |\mu_k(x,\xi,t)| \,d\xi dxdt.
$$
But $\mu_k$ vanishes for $k\to \infty$ because from \eqref{eq:f} we infer that $\mu=0$ both when $f=1$ and $f=0$.
Therefore, we find \eqref{eq:muvanish}.
\qed

\subsection{Proof of Theorem \ref{th:lim}}

The proof of the Theorem \ref{th:lim} can now be easily deduced from 
Lemma \ref{lem:convp}. 
First, up to a subsequence, we have that $p_k$ converges a.e. towards $p_\infty$. 
On the one hand, recalling that the sequence $(p_k)$ is uniformly bounded in $L^\infty$, we use the Lebesgue dominated
convergence Theorem to show that, for any bounded open $U$,
$$
\int_0^T \int_U p_k|p_k-W_k-\nu G(p_k)|\,dx \underset{k\to +\infty}{\longrightarrow}
\int_0^T \int_U p_\infty|p_\infty-W_\infty-\nu G(p_\infty)|\,dx.
$$
On the other hand, we have from estimate \eqref{estim1} that
$$
\int_0^T \int_U p_k|p_k-W_k-\nu G(p_k)|\,dx \underset{k\to +\infty}{\longrightarrow} 0.
$$
We deduce that $p_\infty\big(p_\infty-W_\infty-\nu G(p_\infty)\big)=0$ a.e. that is \eqref{relcomp}.

We may apply the strong convergence for transport equations, as in \cite{Jabin, DiPernaLions}, 
to conclude that, since the term $G(p_k)$ converges strongly, $n_k$, which solves the transport equation
\eqref{eqvis:n}, itself converges strongly.
Note in particular that, from assumption \eqref{id2}, we have 
$n^0_k \underset{k\to +\infty}{\longrightarrow}  \mathbf{1}_{\Omega^0}$.
Passing to the limit in the equation \eqref{eqvis:n}, we recover the limit equation for 
$n_\infty$ \eqref{eq:nlim}.

Finally, passing to the limit in the relation, 
$$
n_k p_k = \Big(\f{k}{k-1}\Big)^{1-\f{1}{k-1}} p_k^{k/(k-1)},
$$
we deduce that $(1-n_\infty) p_\infty=0$.
The relation \eqref{eq:plim} is then a direct consequence of Lemma \ref{lem:convp}.
\qed

\section{One dimensional traveling waves}
\label{sec:tw}

In order to examplify Theorem \ref{th:lim}  and to give a simple case, with  a solution that can be build analytically, we look for a one dimensional traveling wave solution to the Hele-Shaw limit. 

Because, traveling waves are defined up to a translation, we may set, in the moving frame, $\Omega(t)=\R_+$. Then, the system rewrites
\begin{equation}\label{eqtwpos}
 p=0,\qquad
-\sg n' - (n W')' = n G(0), \qquad - \nu W'' + W = 0, \quad \text{for }  x>0, 
\end{equation}
\begin{equation}\label{eqtwneg}
n=1, \qquad -\nu W'' + W - H(W) = 0,
\qquad p = H(W),  \quad \text{for }  x < 0.
\end{equation}
Moreover, the jump condition at the interface $x=0$ implies
$-\sg [n] - [n W'] = 0$, which leads to the traveling velocity
$$
\sg = - W'(0).
$$
We denote $W_0:=W(0)$. For $x>0$,  we have 
\beq\label{Wxpos}
W(x) = W_0 e^{-x/\sqrt{\nu}},
\eeq
from which we deduce that 
$$
\sigma = \frac{W_0}{\sqrt{\nu}}.
$$
Then we can rewrite the first equation in \eqref{eqtwpos} as
$$
-n'(x) \Big(\f{W_0}{\sqrt{\nu}} - \f{W_0}{\sqrt{\nu}} e^{-x/\sqrt{\nu}}\Big)
= n(x)  \Big(G(0)+ \f{W_0}{\nu} e^{-x/\sqrt{\nu}}\Big).
$$
Taking the limit $x\to 0$ leads to $n(0)=0$. Moreover, since $n'\leq 0$, 
we deduce that $n=0$ on $(0,+\infty)$.

For $x<0$, we solve the second order ODE for $W$ with boundary condition
$W(0)=W_0$ and $W'(0)=-W_0/\sqrt{\nu}$.
As an example, we choose for the growth term the function 
\beq\label{Gex}
G(p)=P_M-p, \quad \text{and thus }  H(W)= \frac{W+\nu P_M}{1+\nu}.
\eeq
Then equation \eqref{eqtwneg} for $W$ rewrites~:
$$
-(\nu+1) W'' + W=P_M.
$$
The only bounded solution on $(-\infty,0)$ such that $W(0)=W_0$ is given by
$$
W(x) = P_M + (W_0-P_M) e^{x/\sqrt{\nu+1}}.
$$
Moreover, the continuity of the derivative implies, from \eqref{Wxpos}, that 
$W'(0)=-W_0/\sqrt{\nu}$. We deduce the value for $W_0$~:
$$
W_0 = \frac{\sqrt{\nu}}{\sqrt{\nu}+\sqrt{\nu+1}}P_M.
$$
Then we conclude that for $x<0$,
$$
W(x)=P_M\Big( 1-\f{1}{1+\sqrt{\nu/(\nu+1)}} e^{x/\sqrt{\nu+1}}\Big).
$$
The pressure is then given by~:
$$
p(x) = P_M\Big( 1-\f{1}{\nu+1+\sqrt{\nu(\nu+1)}} e^{x/\sqrt{\nu+1}}\Big).
$$
and the traveling velocity
$$
\sigma = \f{P_M}{\sqrt{\nu}+\sqrt{\nu+1}}.
$$
We notice that the pressure is nonnegative and has a jump at the interface $x=0$.
The height of the jump is given by $P_M\big(1-1/(\nu+1+\sqrt{\nu(\nu+1)})\big)$.
We observe moreover that $\sigma$ is a decreasing function of $\nu$.
Letting $\nu \to 0$, we recover the result for the 
Hele-Shaw model for purely elastic tumors \cite{PTV_TWHS, TVCVDP}.
\\
\\

\begin{figure}[!ht]
\begin{center}
\includegraphics[width=8cm]{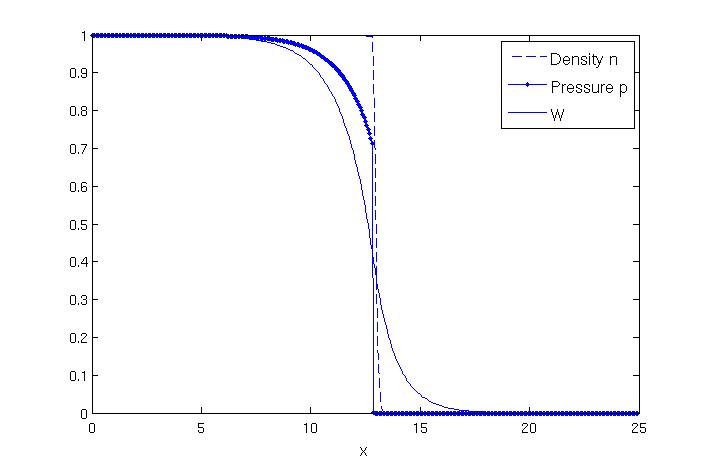}
\includegraphics[width=8cm]{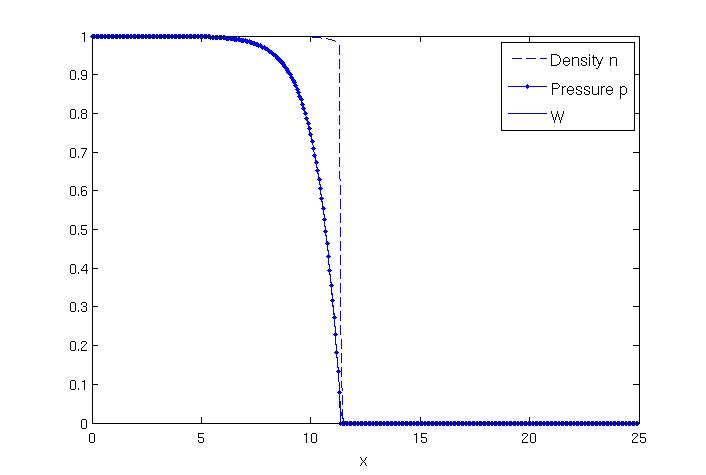}
\end{center}
\vspace{-8mm}
\caption{Plot of the density $n$ (dashed line), pressure $p$ (line with dot), $W$ (continuous line).
Left : for $\nu=1$ and at final time $T=25\,s$. 
We notice a jump for the density from 0 to 1 at the front and a jump of the pressure.
Right : for $\nu=0$ and at final time $T=12.5\,s$. In this case, we have $p=W$ and there is no jump on the
pressure; moreover, the velocity of the front is faster. 
This observation is compatible with the interpretation that viscosity acts as a friction.}
\label{fig:density}
\end{figure}

\begin{figure}[!ht]
\begin{center}
a) \includegraphics[width=7cm]{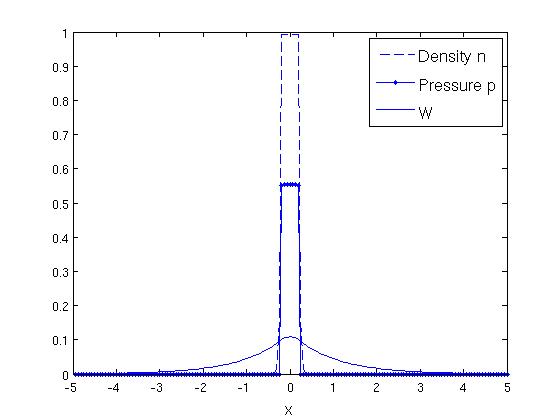}
b) \includegraphics[width=7cm]{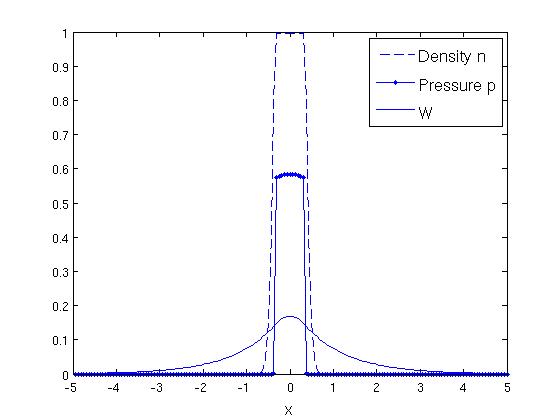} \\
c) \includegraphics[width=7cm]{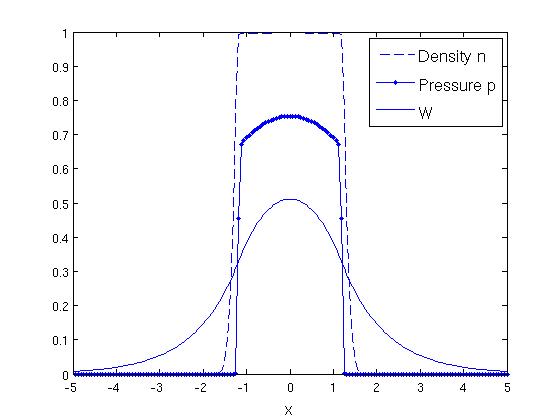}
d) \includegraphics[width=7cm]{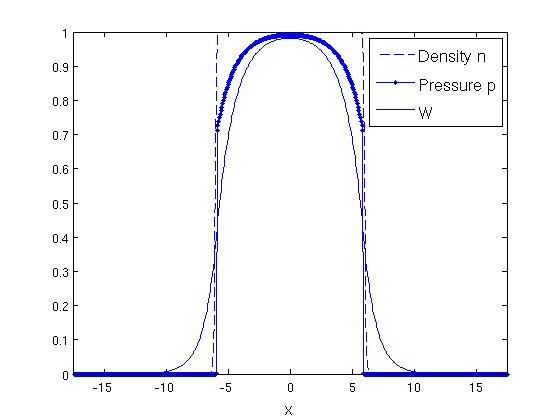}
\end{center}
\vspace{-8mm}
\caption{First steps of the formation of the propagating front with $k=100$ and $\nu=1$.
The density $n$ (dashed line), the pressure $p$ (line with dot) and the potential $W$ (continuous line) 
are represented at 4 successives times : 
a) $t=0.1\,s$, b) $t=1.25\,s$, c) $t=3.75\,s$ and d) $t=12.5\,s$.}
\label{fig:dyna}
\end{figure}

{\bf Numerical simulations.} 
Finally, we present numerical simulations of the system 
\eqref{eqvis:n}--\eqref{pneg} in one dimension.
We use a discretization thanks to a cartesian grid of a bounded domain $[-L,L]$ of the real line.
Equation \eqref{eqvis:n} is discretized by a finite volume upwind scheme. 
Equation \eqref{eqvis:W} is discretized thanks to finite difference scheme.
Since we focus on the case where $k$ is large, we use $k=100$ in the numerical computation.
For the initial data, we choose $n^0=\mathbf{1}_{[-0.2,0.2]}$.
The growth function $G$ is chosen as in \eqref{Gex} with $P_M=1$. 

In Figure \ref{fig:density}, we display the shape of the density $n$, the pressure $p$ and $W$ 
obtained by the numerical simulation. The figure on the left displays the result 
with a viscosity coefficient $\nu=1$. For the comparison, we plot on the right
of Figure \ref{fig:density}, the shape in the case without viscosity ($\nu=0$).
Comparing both figures, we observe that in the case $\nu=1$, we have a jump of the
pressure at the interface of the solid tumor, whereas in the case $\nu=0$,
the pressure is continuous at the interface.

We display in Figure \ref{fig:dyna} the first steps of the formation of the propagating front
with the initial data $n^0=\mathbf{1}_{[-0.2,0.2]}$. For this simulation we take $\nu=1$ and $k=100$.
The dynamics is represented thanks to the plot at 4 successives times of the density $n$,
pressure $p$ and $W$.
After a transitory regime during which the pressure increases until reaching its maximal value $P_M=1$,
the shape of the traveling waves is obtained and the front of the tumor invades the healthy tissue.

\section{Conclusion}

A geometric model, also called incompressible,  has been derived from a cell density model (also called compressible) when the pressure law is stiff.  Because the viscosity is considered here, the limiting problem is a free boundary problem for the set $\Omega(t)$ of non-zero pressure. The limiting system for the pressure consists in an algebraic 
relation between the pressure $p_\infty$ and the limiting potential $W_\infty$ \eqref{relcomp},
coupled with an elliptic equation for the potential $W_\infty$ set in the whole space \eqref{eq:Wlim}.

This is a major difference with the case where viscosity is neglected, the so-called Hele-Shaw system \cite{friedman, Lowengrub_survey}; then, the pressure is given  by an elliptic equation for the pressure in the  moving domain  $\Omega(t)$. A paradox is that the effect of keeping viscosity generates a jump of the pressure at the interface of the region defining the tumor, unlike in \cite{PQV,PQTV} where the Hele-Shaw problem is complemented
with Dirichlet boundary conditions and therefore the pressure is continuous.
This point is also observed in the numerical simulations in Section \ref{sec:tw}.
The velocity of the propagating front of the tumor is given by the equation satisfied
by the density \eqref{eq:nlim}. 
Because  the pressure is discontinuous, it has weaker regularity that in the inviscid case treated in \cite{PQV,PQTV}
and we need to develop a new strategy of  proof to derive  the incompressible limit. Our approach is based on a kinetic formulation of the equation satisfied by the pressure.
\\

This work also opens several additional questions. First, the case of general initial
data is not treated here because we assume that $n^0$ vanishes outside $\Omega^0$. Then, it would be interesting to consider the case with active 
motion as in \cite{PQTV}. In such a case, equation \eqref{eqvis:n} is replaced by 
a parabolic equation. Then the structure of the problem is different but the limiting system
should be the same, except the equation for the density which implies then a faster
propagation of the region $\Omega(t)$.
Finally, it is formally clear from \eqref{omega:nlim}--\eqref{omega:plim} that
letting $\nu\to 0$, we recover the Hele-Shaw system. However, a rigorous proof of this
fact requires compactness of the sequence which we is not directly available with the method developed here.

\bigskip
\noindent {\em Acknowledgment.} This work has been supported by the french "ANR blanche" project Kibord:  ANR-13-BS01-0004.

%
%
%

\end{document}